\documentclass[a4paper,10pt]{article}

\usepackage{float} 
\usepackage{floatflt} 
\usepackage{array}

\usepackage{amssymb}
\usepackage{amsmath}
\usepackage{theorem}
\usepackage[french]{babel}     
\usepackage[ansinew]{inputenc} 

\renewcommand{\emptyset}{\varnothing}
\renewcommand{\epsilon}{\varepsilon}
\renewcommand{\Re}{\operatorname{Re}}
\renewcommand{\Im}{\operatorname{Im}}

\newcommand{\on}{\operatorname}

\newcommand{\ds}{\displaystyle}

\newcommand{\dans}{\rightarrow}
\newcommand{\tend}{\longrightarrow}

\newcommand{\ic}{\mathbf{i}}

\newcommand{\bbC}{\mathbb{C}}
\newcommand{\bbD}{\mathbb{D}}

\newcommand{\bbN}{\mathbb{N}}

\newcommand{\bbR}{\mathbb{R}}
\newcommand{\bbS}{\mathbb{S}}

\newcommand{\bbU}{\mathbb{U}}

\newcommand{\bbZ}{\mathbb{Z}}

\newcommand{\mcZ}{\mathcal{Z}}

\usepackage{a4}        
\usepackage{ifthen}
\usepackage{calc}      
\usepackage{epic}      
\usepackage{eepic}     
\usepackage{graphicx}  
\usepackage{rotating}  
\usepackage{float}     
\usepackage{rotating}  
\usepackage{amscd}     
\usepackage[all]{xy}   
\usepackage{xspace}    

\newtheorem{defi}{Définition} 
\newtheorem{prop}[defi]{Proposition}

\newtheorem{lem}[defi]{Lemme}

\newenvironment{dem}{\noindent{\bf Preuve}~:}%
                    {\hfill{$\blacksquare$}\par\medskip}

{\begin{list}{--}{\setlength{\itemsep}{1.3ex}%
                  \settowidth{\labelwidth}{--}%
                  \setlength{\leftmargin}{\labelwidth+\labelsep-\itemindent}}
}%
{\end{list}}

\renewcommand{\Im}{\on{Im}}
\renewcommand{\Re}{\on{Re}}
\renewcommand{\epsilon}{\varepsilon}

\begin{document}

\setlength{\parindent}{0pt}

\begin{center}
  \huge
  Cage de Faraday
\end{center}

\bigskip

\begin{center} Version du 5 août 2002. \end{center}

{\tiny Avertissement~: cet article n'a pas un statut fixe et est
  susceptible d'évoluer à tout moment.}

\vspace{2.5cm}

\textbf{\LARGE I. Introduction}

\bigskip

\bigskip

\textbf{\Large Motivation}

\bigskip

Dans ma thèse~\cite{C}, la preuve du lemme~5.4 est fausse (page~40,
partie~I, section~5~: ``Application à une nouvelle preuve du théorème
de Yoccoz''). Cela remet en question le lemme~5.5 et son
corollaire~5.6 qui suivent.

Je donne ici un énoncé analogue au corollaire~5.6 qui suffit à
l'utilisation qui en est faite dans ma thèse, ainsi qu'une preuve
inspirée du travail de Rohde et Zinsmeister~\cite{RZ}.

Notons que Xavier~Buff a trouvé une preuve très élégante d'un énoncé
plus fort que le lemme~5.5, ce qui sauve ce dernier ainsi que son
corollaire~5.6.

\bigskip

\bigskip

\textbf{\Large Situation du problème}

\bigskip

Dans~\cite{C}, nous composons une suite de revêtements universels $\phi_n
: (\bbD,0) \dans (\bbD\setminus X_n,0)$ où $X_n$ est un ensemble fini
de points non nuls. La dérivée en $0$ de la composée $\phi_n \circ \cdots
\circ \phi_0$ est de valeur absolue décroissante en fonction de $0$, et
nous cherchons à prouver qu'elle tend vers $0$ quand $n \tend
+\infty$. Nous disposons pour cela d'informations sur la suite $X_n$~:
$X_n$ est l'image par une certaine fonction $f_n : \bbD \dans \bbD$
fixant $0$ de l'ensemble $r_n \bbU_{q_n}$ où $r_n \in ]0,1[$, $q_n
\in \bbN^*$, et $\bbU_q$ désigne l'ensemble des racines $q$-ièmes de
l'unité. De plus, nous savons que $\prod r_n \tend 0$ et que la
distance hyperbolique dans $\bbD$ de $r_n$ à $r_n \exp(\ic 2\pi /q_n)$
tend vers $0$. Ces conditions sont-elles suffisantes~?

\

Dans cet addendum, nous ne retiendrons que le fait qu'il
existe une constante $d>0$ tel que pour tout $\bbN$, $X_n$ est une
chaîne cyclique de points de $\bbD$, chacun situé à distance distance
hyperbolique $\leq d$ du suivant, faisant le tour de $0$ au moins une
fois, et tous à distance euclidienne $\geq 1-r_n$ du bord de $\bbD$.
Nous obtenons alors une majoration du type
\[\ln |\phi'(0)|\ <\ \frac{\ln r_n}{K(d)}\ \ (< 0)\]
pour une certaine constante $K(d)>0$ dépendant de $d$.

\

Note~: dans le même temps que le présent travail, Xavier Buff à
démontré que sous la condition originelle $X_n = f(r_n \bbU_{q_n})$,
alors
\[|\phi_n'(0)| < |\psi_n'(0)|\]
où $\psi_n$ est le revêtement universel $(\bbD,0) \dans (\bbD
\setminus r_n\bbU_{q_n} , 0)$, or il est relativement facile d'obtenir une
majoration du type $\ds \ln |\psi_n'(0)| < \frac{\ln r_n}{K'(d)}$ (voir
les calculs effectués dans~\cite{C}, partie~I, section~6.5, pages~60
et~61, montrant une majoration), ce qui est suffisant. De plus, $K'(d)
\tend 1$ quand $d\tend 0$.

\bigskip

\bigskip

\textbf{\Large Heuristique}

\bigskip

Prenons le disque $\bbD$ et retirons-lui un seul point~: $X=\{z\}$. On
peut calculer explicitement le rayon conforme $r$ (voir~\cite{C}) de
$\bbD \setminus X$. Notons $\epsilon = 1-|z|$, et faisons tendre $|z|$
vers $1$. Alors $r = 1 -\epsilon^2/6 + O(\epsilon^3)$. Ainsi la
contribution d'un point proche du bord du disque est de l'ordre du
carré de la distance au bord. Si maintenant, on prend pour $X$ la
composante contenant $0$ du disque privé de la géodésique passant par
$z$ et orthogonale à $[0z]$, alors $r = 1 -\epsilon^2/2 +
O(\epsilon^3)$. Ainsi la chute de rayon conforme est du même ordre de
grandeur que pour le point isolé. Dans le second cas, c'est également
l'ordre de grandeur de l'aire de la partie retirée. Notons
momentanément $H(z)$ la région retirée.

Si on retire une courronne concentrique, la perte est également de
l'ordre de l'aire retirée. Si on retire une suite finie de points
régulièrement répartie sur un cercle concentrique, tels que la
distance entre deux points consécutifs est de l'ordre de la distance
au bord, la perte de rayon conforme sera de l'ordre de cette distance
(voir~\cite{C}, partie~I, section~6.5), donc de l'ordre de l'aire de
la réunion des régions $P(z)$.

Si on retire un nuage de points très proches d'un même point
$z\in\bbD$, alors d'une part le rayon conforme est quasiment le même
qu'en retirant le seul point, d'autre part les régions $P(z)$ se
superposent, et donc à nouveau l'aire de leur réunion et la perte de
rayon conforme sont commensurables.

Est-ce vrai plus généralement~? C'est l'objet de la
proposition~\ref{prop_CRZ} de cet article. Notons que pour des raisons
pratiques, la forme des régions retirées y diffère.

\clearpage

\textbf{\LARGE II. Rayon conforme des grands ouverts du disque}

\bigskip

\medskip

\begin{figure}[htbp]
 \begin{center}
  \begin{picture}(280,195)
   \scalebox{0.5}{
    \put(0,0){\includegraphics[300pt,200pt][410pt,590pt]{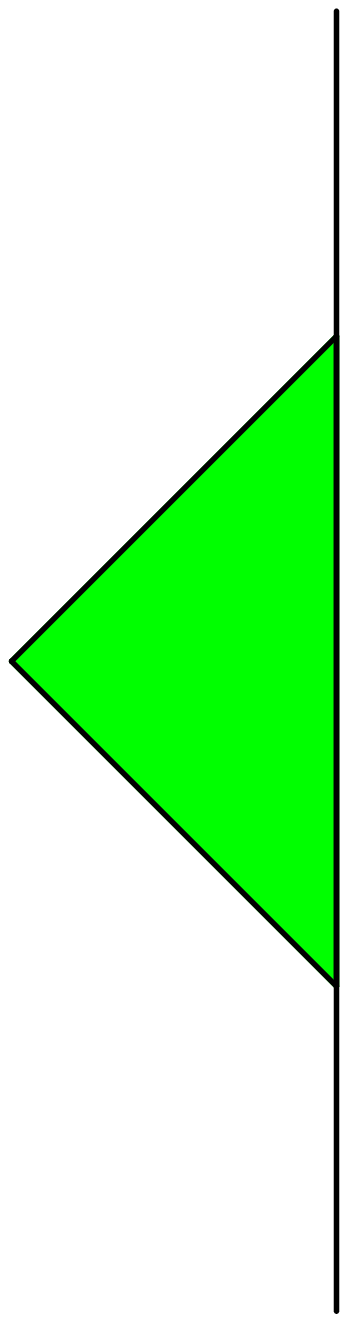}}
    \put(190,0){\includegraphics[120pt,200pt][490pt,590pt]{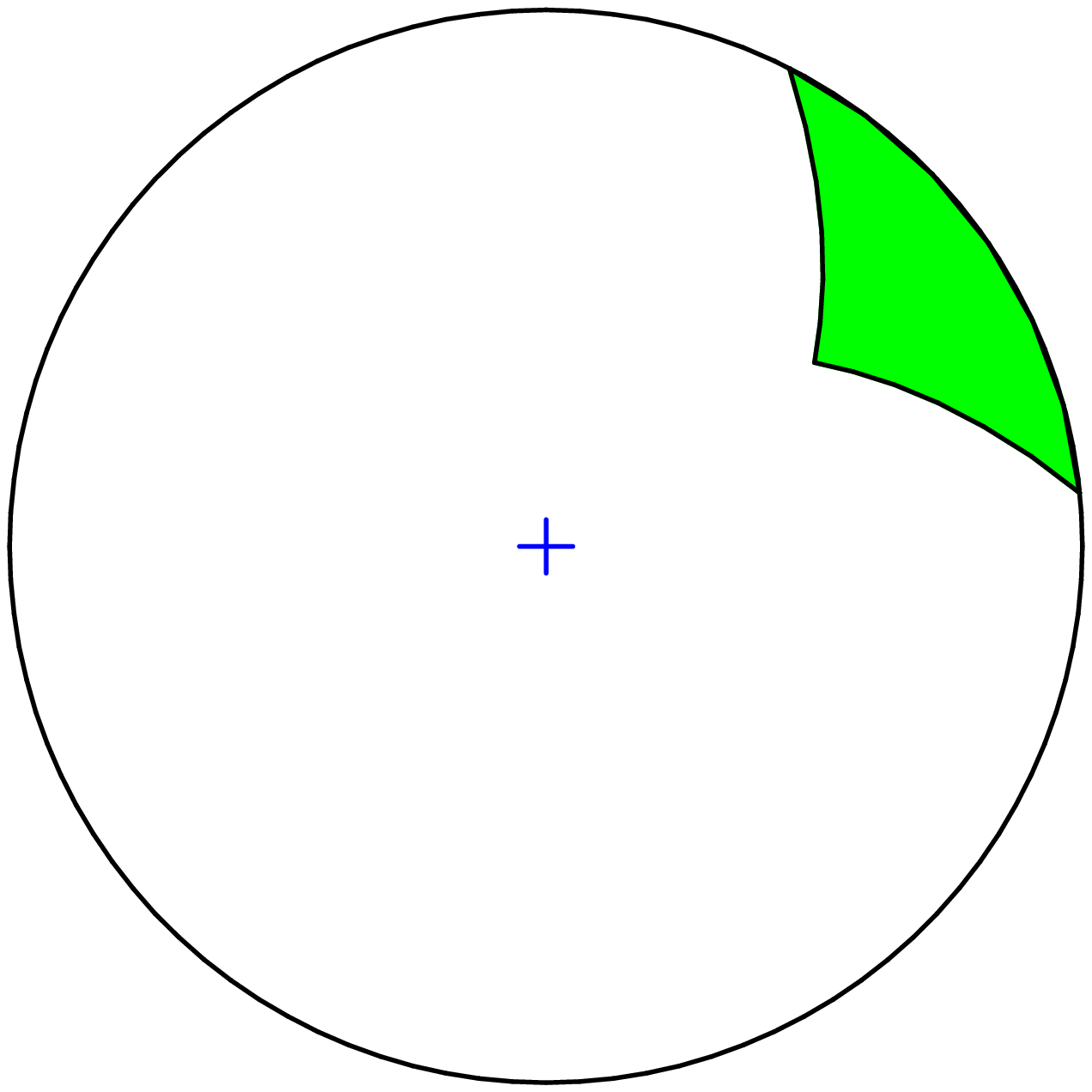}}
   }
   \put(-10,98){$w_0$}
   \put(30,98){$T$}
   \put(70,98){$\overset{\exp}{\longrightarrow}$}
   \put(220,125){$z_0$}
   \put(248,140){$V$}
  \end{picture}
 \end{center}
\end{figure}

Pour $z_0 \in \bbD$ non nul,
soit $w_0$ une quelconque préimage par l'application exponentielle, et
considérons le triangle $T$, intersection du demi-plan ``$\Im(w)\leq 0$''
et du secteur d'équation $|\on{Arg}(w-w_0)|\leq \pi/4$, sommet compris
(triangle équilatéral isocèle d'hypoténuse posée sur l'axe imaginaire
et de sommet droit placé en $w_0$). Soit $V(z_0)=\exp(T)$.

Notons que le domaine $V(z)$ s'auto-intersecte si $|z|\leq e^{-\pi}$. Ce
ne sera pas un problème, d'autant plus que les points $z$ qui nous
intéresseront seront proches du bord du disque.

\begin{figure}[htbp]
 \begin{center}
  \begin{picture}(310,195)
   \put(0,5){\scalebox{0.333}{
    \includegraphics[200pt,130pt][430pt,660pt]{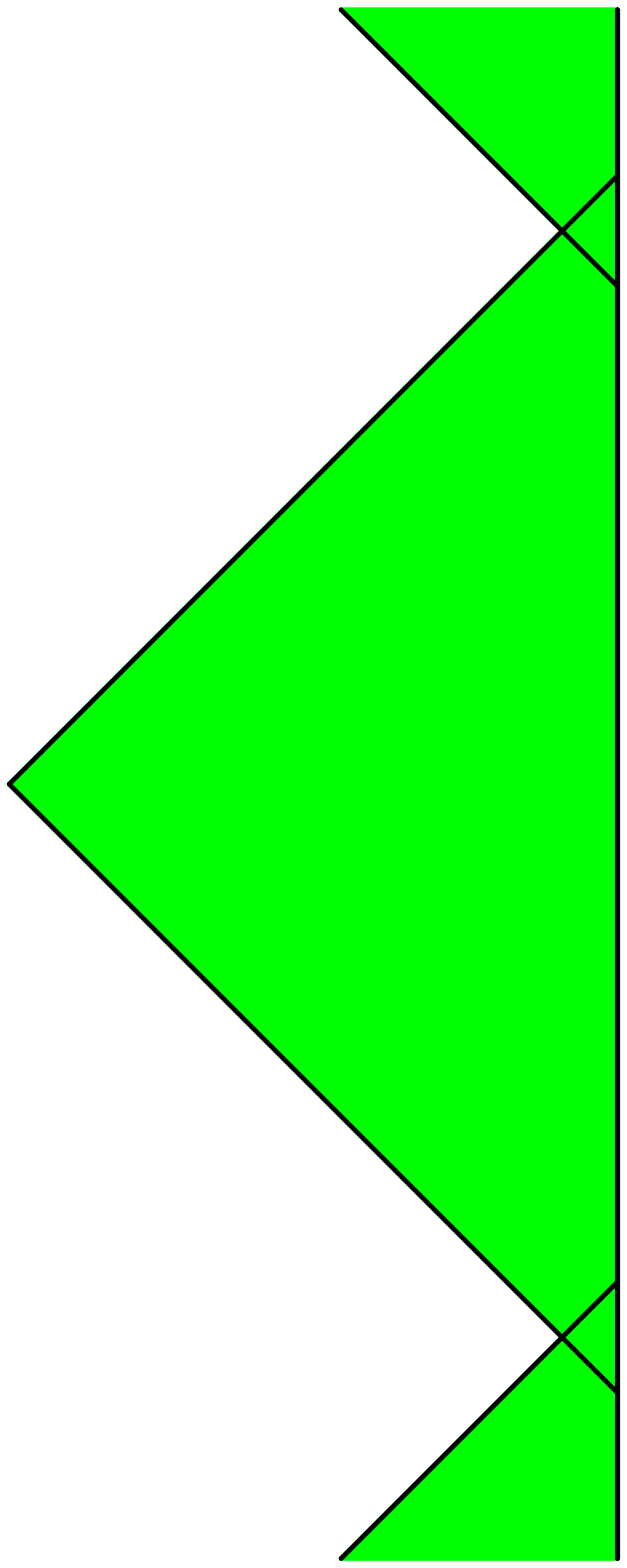}
   }}
   \put(90,95){$\overset{\exp}{\longrightarrow}$}
   \put(115,0){\scalebox{0.5}{
    \includegraphics[120pt,200pt][490pt,590pt]{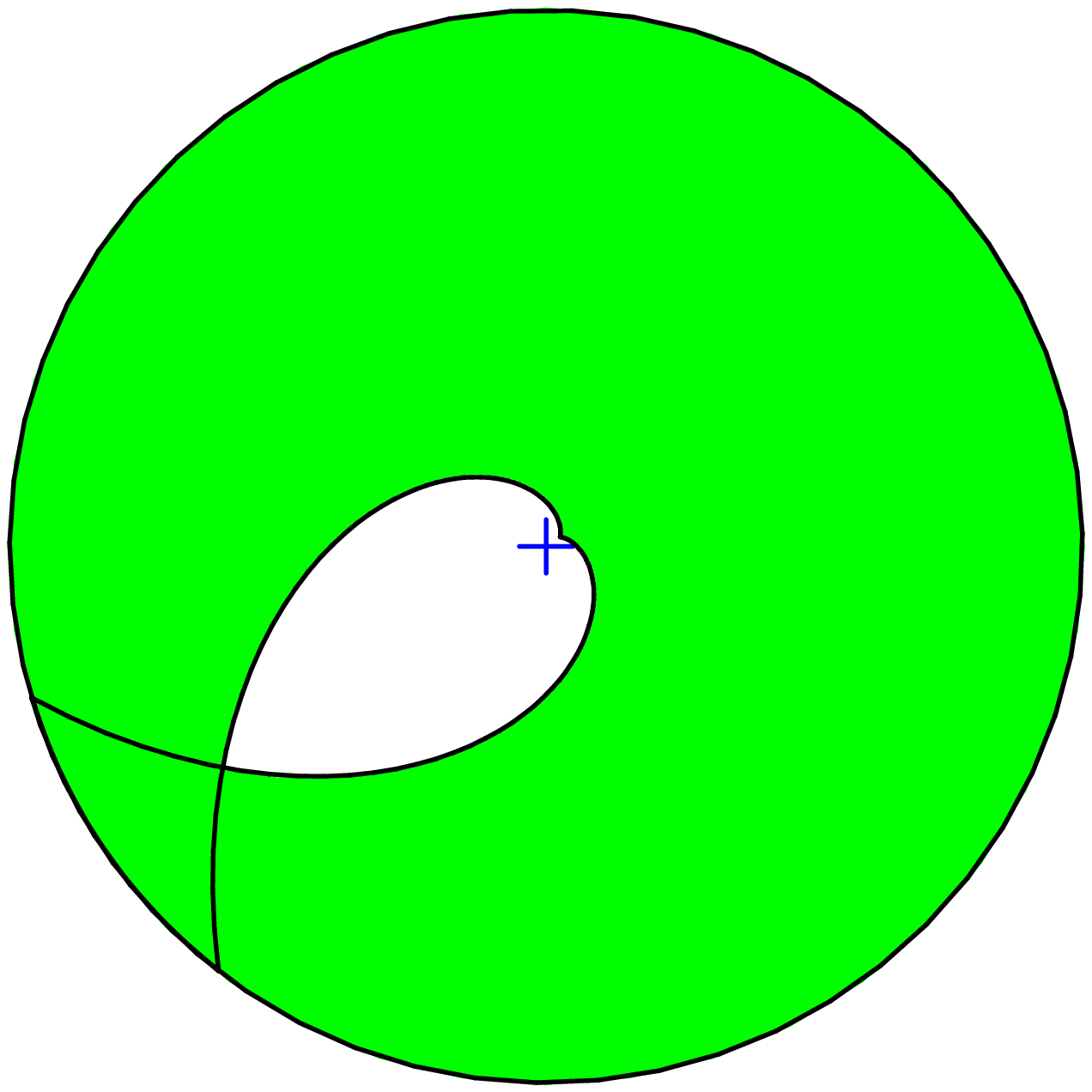}
   }}
  \end{picture}
 \end{center}
\end{figure}


Remarquons également que le choix d'un triangle en coordonnées
logarithmiques est arbitraire, et qu'il y a d'autre formes et conventions
qui conviendraient pour ce qui va suivre.

Avant d'énoncer la proposition~\ref{prop_CRZ}, intéressons nous
d'abord à un cas simple.
\begin{lem}[contribution d'un point]\label{lem_calculs}
  Si $u>0$, $x=e^{-u}$, $U=\bbD\setminus\{x\}$, et $r$ est le rayon
  conforme associé à $U$, alors
  \[r = \frac{u}{\sinh(u)}\]
  Quand $u \tend 0$,
  \[r = 1 -\frac{u^2}{6} + O(u^3)\]
  pour tout $u \in ]0,\pi]$,
  $-\ln(r)/u^2 \in [K_1,1/6]$, avec $K_1>0.13$.
\end{lem}
Voir ma thèse (\cite{C}) pour les deux premières affirmations.
La troisième est élémentaire.

Soit $u>0$, $x=\exp(-u)$, $r$ le rayon conforme de $\bbD \setminus
\{x\}$, et $a=1/2\pi \times$ l'aire de $T(-u)$ \emph{dans le
quotient} $\bbC/ \ic 2\pi \bbZ$. On calcule explicitement que $\ds r =
\frac{u}{\sinh(u)}$, que $a= u^2/2\pi$ si $u\in]0,\pi]$, $a= u -
\pi/2$ si $u \in [\pi, +\infty[$, et une étude de fonction montre
alors que le quotient $\ds\frac{|\ln r|}{a/2\pi}$ est toujours compris
entre $0.75$ et $1.05$.

\begin{figure}[htbp]
 \begin{center}
  \begin{picture}(360,245)   
    \put(0,60){\scalebox{0.45}{%
      \includegraphics[120pt,200pt][490pt,590pt]{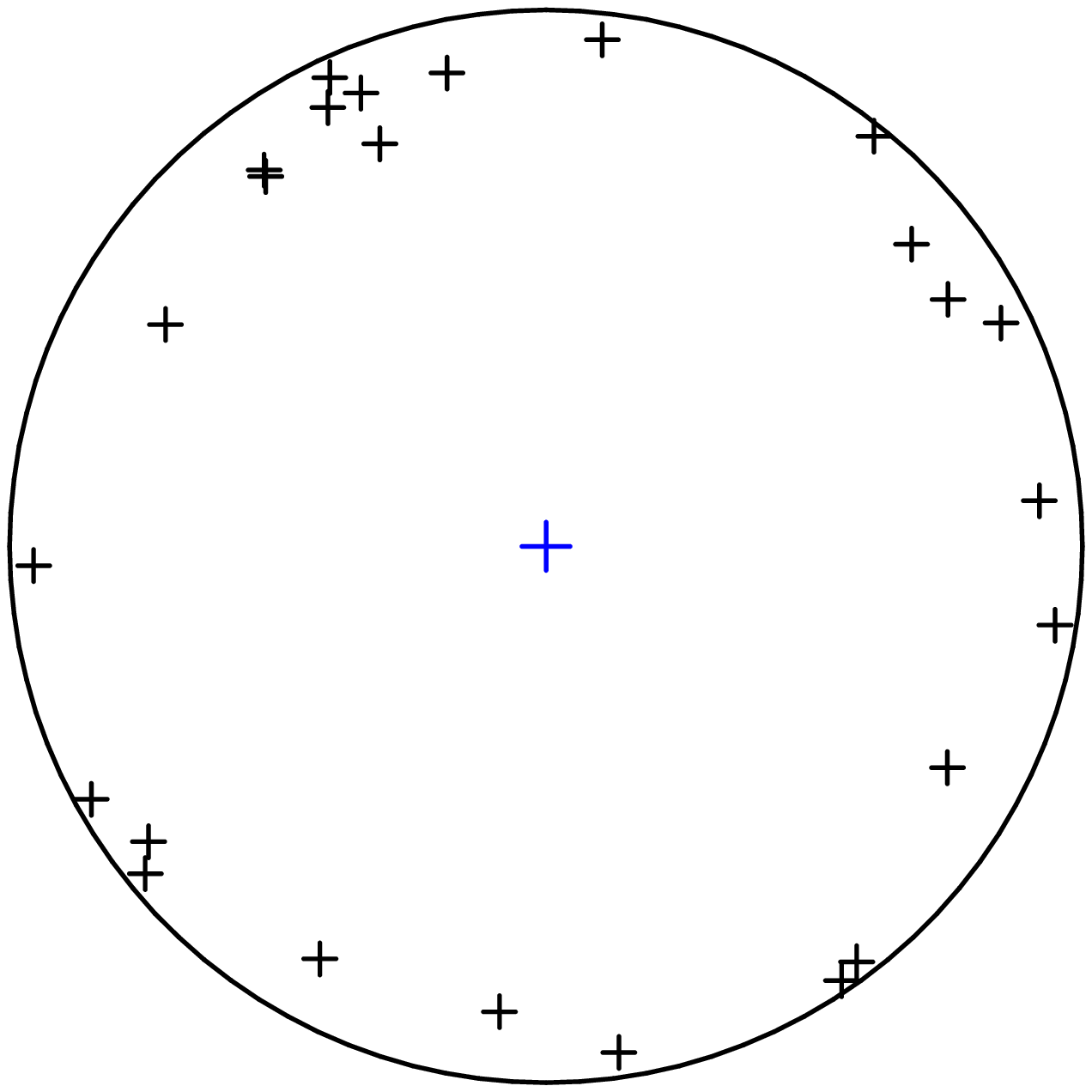}}}
    \put(190,60){\scalebox{0.45}{%
      \includegraphics[120pt,200pt][490pt,590pt]{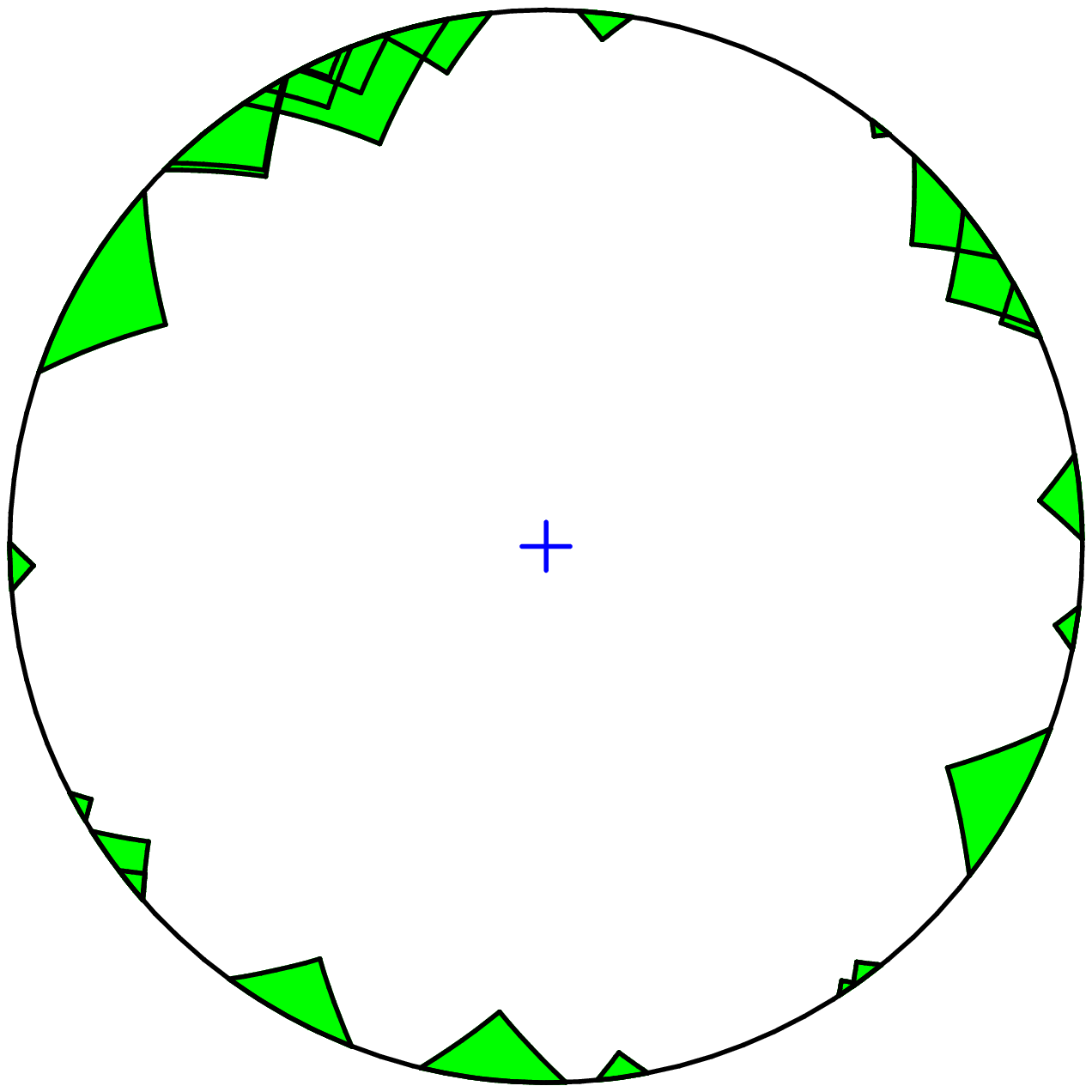}}}
    \put(35,20){\parbox{290pt}{Figure~: illustration de la
    proposition~\ref{prop_CRZ}. Le domaine $\bbD$ privé des
    points est approximé par $\bbD$ privé des ``triangles'' issus de
    ces points. L'aire de la réunion des triangles pleins, divisée par
    $2\pi$, est de l'ordre de $1-r$}}
  \end{picture}
 \end{center}
 \label{fig_ill}
\end{figure}

La proposition suivante généralise un résultat de Rohde et Zinsmeister
(\cite{RZ}).

\begin{prop}\label{prop_CRZ}
  Il existe des constantes $K'>K>0$ telles que pour tout
  ouvert $U$ inclus dans $\bbD$ et contenant $0$, notons
  $r$ le rayon conforme en $0$ du revêtement universel de $U$,
  $\ds A = \bigcup_{e^w \in \bbD \setminus U} T(w)$ et $a =
  \on{aire}(A)/2\pi$ mesurée dans le quotient $\bbC /\ic 2\pi \bbZ$.
  Alors
  \[ -K'.a \leq \ln r \leq -K.a\]
\end{prop}

De façon équivalente, (si $U\not=\bbD$)
\[\ds K \leq \frac{|\ln r|}{a} \leq K'\]


Seule la majoration de $r$ sera utilisée pour la correction de la
thèse.

\

Notons que ce lemme a déja de l'intérêt dans le cas où le domaine $U$
est simplement connexe. Alors, $r$ est son rayon coforme au sens
classique.

\bigskip

\bigskip

\textbf{\large Un lemme utile}

\bigskip

Le lemme suivant sera utilisé deux fois

\begin{lem}\label{lem_vitally_zero}
  Pour tout ouvert
  $O$ inclus dans le cylindre $L=\bbC/2\ic\pi\bbZ$ et
  contenant le demi-plan $``\Re(z) < -\pi"$, il
  existe un ensemble fini ou dénombrable de 
  points $w_i\in O$, $i \in I$, tels que les domaines $T(w_i)$
  soient deux à deux disjoints, et tels que
  \[\ds \bigcup_{w\in L \setminus O} T(w) \subset 
  \bigcup_{i \in I} T_5(w_i)\]
  où $T_5(w) = T(5\Re(w)+\ic\Im(w))$ est un triangle de base $5$ fois
  plus grande que $T(w)$.
  Par conséquent,
  \[\ds \sum_{i \in I} \on{aire}\ T(w_i) \ \leq \ \on{aire}
    \bigcup_{w\in L \setminus O} T(w) \ \leq \ C_3 . \sum_{i \in
    I} \on{aire}\ T(w_i)\] 
  où $C_3 = 25$.
\end{lem}
\begin{dem}
  Construisons par récurrence une suite $W_n$ de sous-ensembles finis
  de points de $L \setminus O$ dont la partie réelle appartient à
  l'intervalle $[-\pi/2^n,-\pi/2^{n+1}[$ de la façon suivante.
  Étant donnés $W_0$, \ldots, $W_{n-1}$ tels que les triangles $T(w)$
  pour $w\in W_0 \cup \cdots \cup W_{n-1}$ sont deux à deux disjoints,
  considérons l'ensemble des points de $L\setminus O$ dont la partie
  réelle appartient à l'intervalle $[-\pi/2^n,-\pi/2^{n+1}[$. Un
  sous-ensemble $W$ tel que les triangles $T(w)$ sont deux à deux
  disjoints possède au plus $2^{n+1}$ éléments, car les segments de
  l'axe imaginaire correspondant à leurs bases doivent être deux à
  deux disjoints. Soit $W_n$ un sous-ensemble maximal tel que les
  triangles $T(w)$ pour $w\in W_n$ sont disjoints des $T(w)$ pour $w
  \in W_0 \cup \cdots \cup W_{n-1}$~: $W_n$ est fini.
  \\ Maintenant, pour tout $w\in L\setminus O$, soit
  $n$ tel que $\Re(w) \in [\pi/2^n,\pi/2^{n+1}[$~: il existe $w' \in
  W_0 \cup \cdots \cup W_n$ tel que $T(w) \cap T(w') \not= \emptyset$.
  La hauteur de $T(w')$ est $>\pi/2^{n+1}$ et celle de $T(w)$ est
  $\leq\pi/2^n$, $T(w) \subset 5T(w')$, où $5T(w')$ désigne le
  triangle homothétique à $T(w')$, de rapport $5$ et de centre le
  milieu de la base de $T(w')$.  On en déduit que l'aire dans le
  cylindre de la réunion des $T(w)$ pour $\ds w\in \bigcup_{n\in\bbN}
  W_n$ est $\leq 25$ fois celle de la réunion des $T(w)$ pour $w\in L
  \setminus O$.
  %
\end{dem}

La construction ``à la Vitalli-Carathéodory'' de cette preuve remplace
l'argument de décomposition dyadique de~\cite{RZ}. En raisonnant plus
finement, on peut faire descendre le facteur $25$ à $5+\epsilon$
pour tout $\epsilon$, et peut-être même en deçà.

\

Voici une version adaptée pour le disque $\bbD$, allégée en vue de son
application à la majoration de $r$~:

\begin{lem}\label{lem_vitalmajor}
  Pour tout ouvert $U$ inclus dans $\bbD$ et contenant le disque
  $B(0,e^{-\pi})$, il existe un ensemble fini ou dénombrable de 
  points $w_i\in U$, $i \in I$, tels que les domaines $T(w_i)$
  soient deux à deux disjoints dans le quotient $\bbC/\ic 2 \pi\bbZ$,
  et tels que
  \[\ds \on{aire} \bigcup_{e^w in\bbD \setminus U} T(w) \leq C_3 .
  \sum_{i \in I} \on{aire} T(w_i)\]
  les aires étant mesurées dans le quotient.
\end{lem}

\

Version allégée pour la minoration de $r$~:
\begin{lem}\label{lem_vitalminor}
  \[\ds \sum_{i \in I} \on{aire} T_5(w_i) \leq  
    C_3 . \on{aire} \bigcup_{e^w\in\bbD \setminus U} T(w)\]
  les aires étant mesurées dans le quotient $\bbC/\ic 2\pi\bbZ$.
\end{lem}

\bigskip

\textbf{\Large Preuve de la proposition~\ref{prop_CRZ}}

\bigskip

Traitons tout de suite le cas où il existe $z \in \bbD\setminus U$ tel
que $|z| \leq e^{-\pi}$. Soit $z$ de module minimal, et $u>0$ tel que
$z= \exp(-u)$. Alors $u - \pi/2 \leq a u$, et $e^{-u}
\ds \leq r \leq \frac{u}{\sinh u}$, d'où
\[\frac{\ln( \sinh(u)/u) }{u}\leq \frac{|\ln r|}{a} \leq
\frac{u}{u-\pi/2}\] 
Une étude de fonction montre alors que
\[0.4 \leq \frac{|\ln r|}{a} \leq 2\]

\

Dans le cas où tous les point $z \in \bbD\setminus U$ vérifient
que $|z| > e^{-\pi}$, nous procéderons en deux étapes, en nous
inspirant de~\cite{RZ}.

\bigskip\bigskip

\textsl{\large Majoration de $r$}

\bigskip

Commençons par quelques calculs.

\begin{floatingfigure}[r]{70pt}
 \begin{picture}(70,160)
  \put(15,0){\scalebox{0.3}{
   \includegraphics[280,130][380,660]{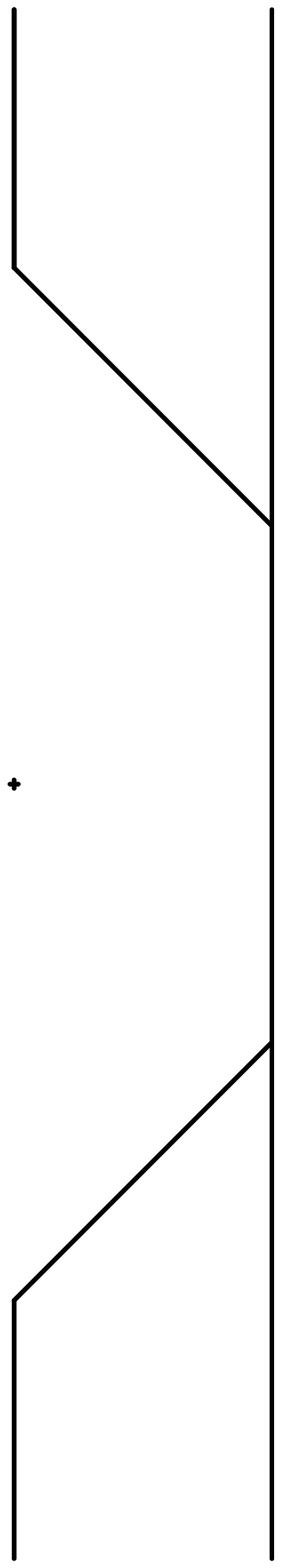}
  }}
  \put(7,78){$w$}
  \put(10,100){$D(w)$}
  \end{picture}
 \label{fig_E}
\end{floatingfigure}

\

\'Etant donné $z=\exp(w)\in\bbD$, soit $E(z)=\exp(D(w))$ où $D(w) =
\Im(w).\ic-\Re(w) 
\times D_0$ et $D_0$ est le sous ensemble de $\bbC$ d'équation $\Re(w)
\leq h(\Im(w))$, avec $h: \bbR \dans \bbR$ une fonction continue
définie par $h(y)=-1$ pour $y\in]-\infty,-2] \cup 
[2,+\infty[$, $h(y)=0$ pour $y\in[-1,1]$, et $h$ affine sur les deux
intervalles restant


\begin{figure}[htbp]
 \begin{center}
  \begin{picture}(375,185)
   \put(-10,0){\scalebox{0.5}{
     \includegraphics[120,210][490,580]{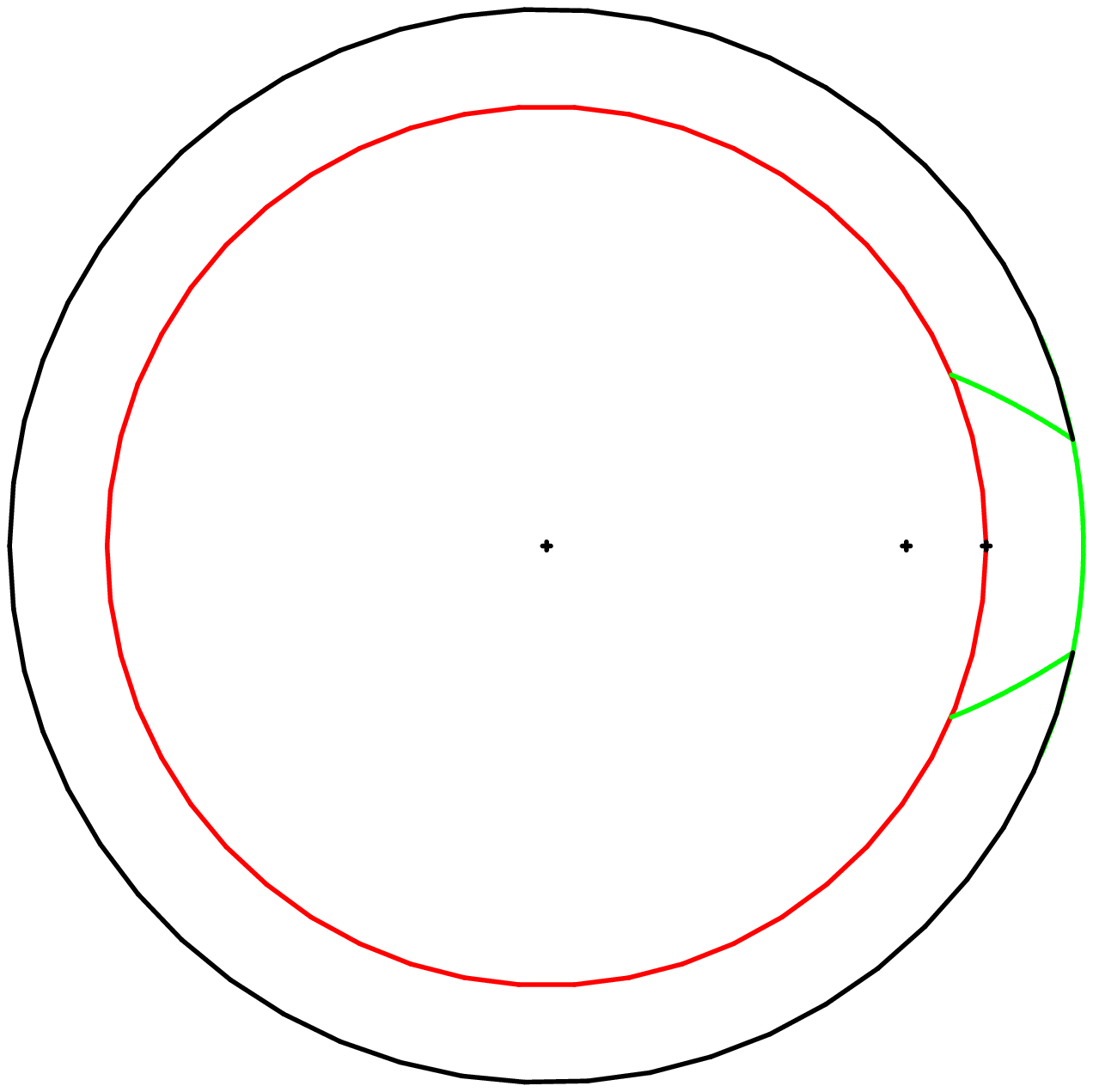}
   }}
   \put(190,0){\scalebox{0.5}{
     \includegraphics[120,210][490,580]{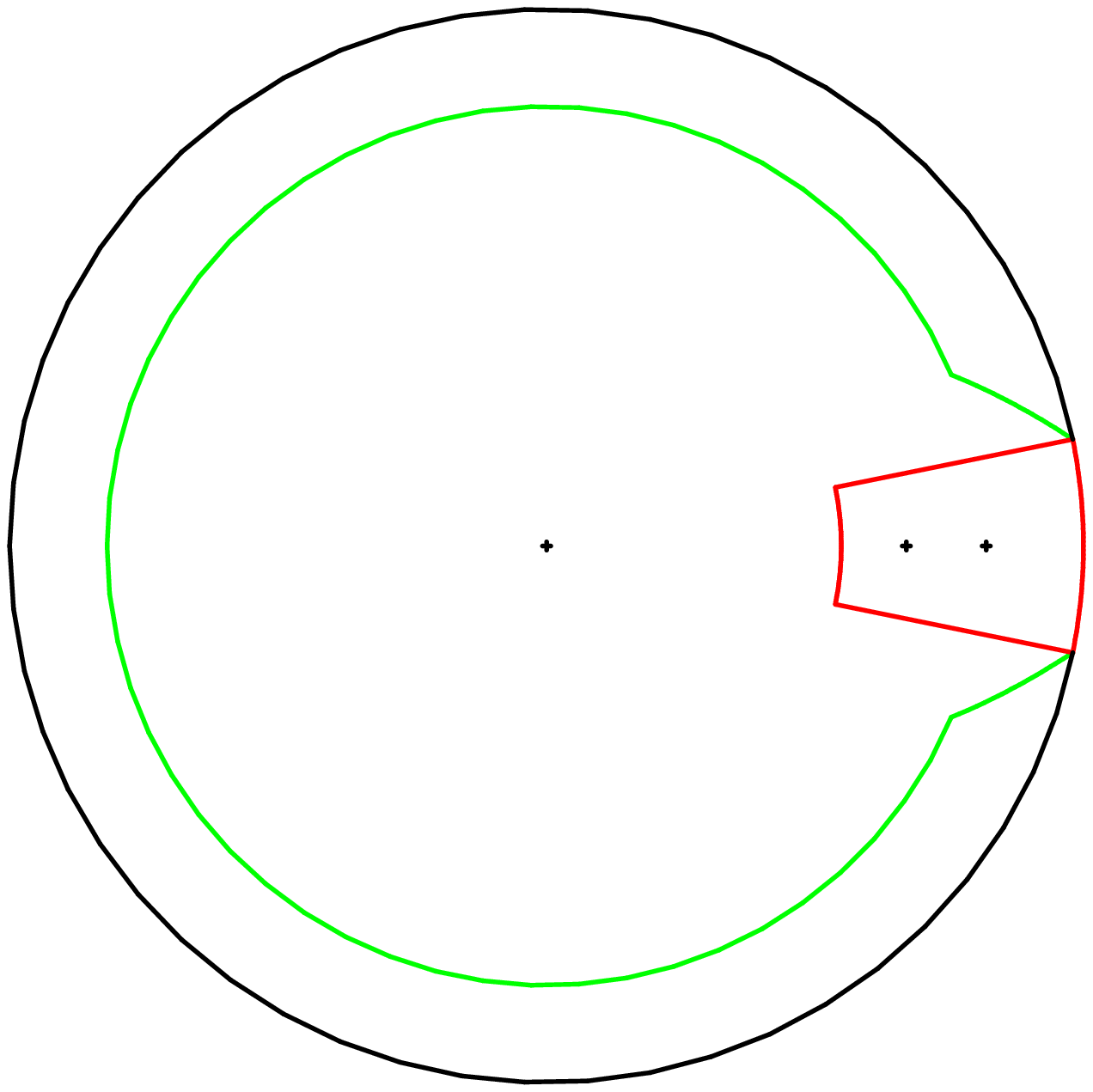}
   }}
   \put(-10,0){
    \put(95,85){$0$}
    \put(153,85){$x^2$}
    \put(168,85){$x$}
   }
   \put(190,0){
    \put(95,85){$0$}
    \put(153,85){$x^2$}
    \put(168,85){$x$}
   }
   \end{picture}
  \label{fig_Es}
 \end{center}
\end{figure}

\begin{lem}\label{lem_E}
  Il existe une constante universelle $C_1$ telle que
  pour tout $x\in[e^{-\pi},1[$, 
  soit $x'$ le point image de $x$ par la représentation conforme de
  $E(x)$ sur $\bbD$ qui envoie $0$ sur $0$ avec dérivée $>0$.
   Alors
  \[\frac{-\ln x'}{-\ln x} > C_1\]
\end{lem}
\begin{dem}
 \begin{floatingfigure}[r]{110pt}
  \begin{picture}(110,215)
   \put(1,0){\scalebox{0.4}{
    \includegraphics[170,130][440,660]{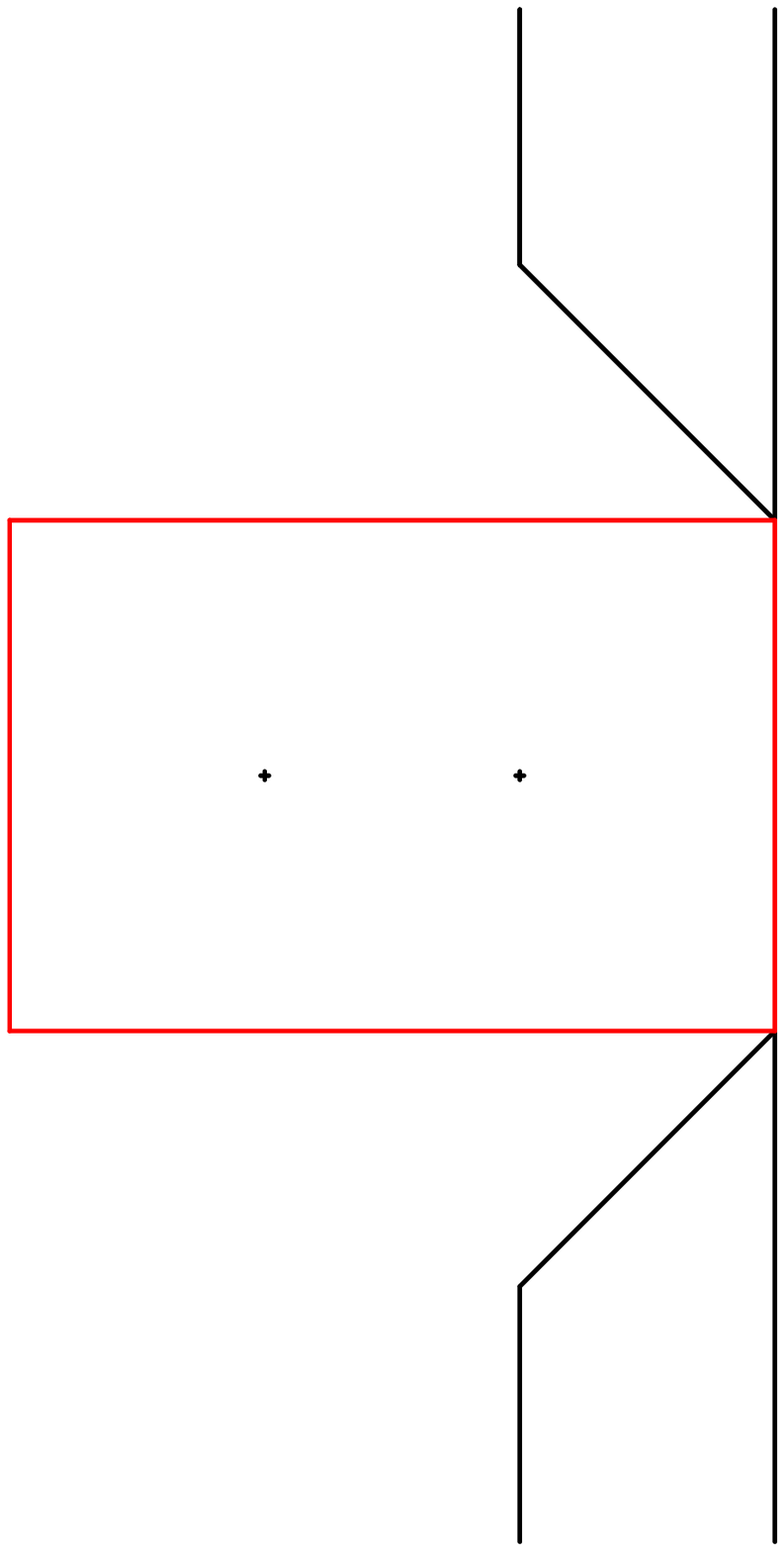}
   }}
   \put(64,92){$\ln(x)$}
   \put(27,92){$\ln(x^2)$}
   \end{picture}
  \label{fig_Rect}
 \end{floatingfigure}
  Majorons la distance hyperbolique $d$ dans $E$ de $0$ à $x$.
  Joignons d'abord $0$ à $x^2$. Comme $E$ contient le disque $B(0,x)$,
  la distance hyperbolique de $0$ à $x^2$ dans $E$ est $\leq$ celle
  dans $B(0,x)$, qui est égale à celle de $0$ à $x$ dans $\bbD$, c'est
  à dire $\tanh^{-1}(x)$. Maintenant, relions $x^2$ à $x$~: $E$
  contient un domaine, image injective par l'exponentielle d'un
  rectangle de même base que le triangle $T(\ln(x))$ et de hauteur
  $3/2$ de celle-ci. Sur ce rectangle, la distance hyperbolique de
  $\ln(x^2)$ à $\ln(x)$ est une constante universelle $C_2$, et il est
  facile de voir que $C_2 < 1$.
  Donc $d \leq \tanh^{-1}(x) + C_2$. Si on uniformise
  $E$ sur $\bbD$, le point $x$ est envoyé sur un point $x'$ situé à
  distance hyperbolique de $0$ égale à $d$, donc à distance
  euclidienne égale à $\tanh(d)$. On en déduit par une suite de calculs
  élémentaires que $\ds \frac{-\ln x'}{-\ln x} \geq
  \exp(-2 C_2) \frac{1-e^{-\pi}}{\pi} > 0.04$ (nous n'avons pas été
  très fins
  dans cette minoration).
\end{dem}

Nous considérons maintenant une suite décroissante $U_i$ de
sous-ensembles de $\bbD$.
Pour simplifier l'exposé, nous supposerons l'ensemble des points 
des points $z_i$ fini. Nous voulons retirer au disque un par un les
triangles associés aux $z_i$ en commençant par les plus petits
(dans le cas où l'ensemble des $z_i$ est dénombrable, leur norme tend
vers $1$ et donc il faut procéder dans l'ordre inverse~: rajouter un à
un les triangles au domaine privé de tous les triangles)
Classons les points $z_i$ de sorte que leur normes forment
une suite décroissante. Soit
$U_0=\bbD$ et $U_{i} 
= U_{i-1} \setminus \{z_{i}\}$, c'est à dire que l'on retire
successivement à $\bbD$ les 
points en commençant par les plus proches du bord. Soit $r_i$ le rayon
conforme du revêtement universel de $U_i$~: $r_0=1$. 
Pour tout $i>0$, le domaine $U_{i-1}$ contient $E(z_i)$. Soit $\phi :
(\bbD,0) \dans (U_{i-1},0)$ revêtement universel. Soit $\psi$ la
branche inverse définie sur $E(z_i)$ et envoyant $0$ sur $0$. Soit
$x''_i = |\psi(z_i)|$. Alors $x''_i<x'_i$ où $x'_i$ est le nombre
correspondant à $x=|z_i|$ dans le lemme~\ref{lem_E}.
Donc, d'après les lemmes~\ref{lem_E} et~\ref{lem_calculs},

\begin{figure}[htbp]
 \begin{center}
  \begin{picture}(185,185)
   \put(0,0){\scalebox{0.5}{
     \includegraphics[120,210][490,580]{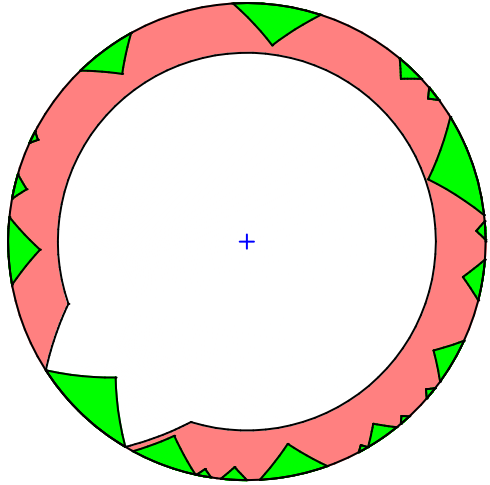}
   }}
   \end{picture}
   \label{fig_superpose}
 \end{center}
 \caption{Les triangles précédent le triangle venant d'être posé sont
     tous en dehors de}
\end{figure}

{
\setlength{\extrarowheight}{8pt}
\[\begin{array}{rclcl}
 \ds - \ln \frac{r_{i}}{r_{i-1}}
 & \geq & \ds K_1 (\ln x''_i)^2 & \quad & 
   \text{(lemme~\ref{lem_calculs})} \\
 & \geq & \ds K_1 (\ln x'_i)^2 & \quad &
   \\
 & \geq & \ds C_1^2 K_1 (\ln |z_i|)^2 & \quad &
   \text{(lemme~\ref{lem_E})} \\
 & = & \ds C_1^2 K_1 \on{aire} T(w_i) & \quad &
    \\
\end{array}\]
}

Soit $K_2 = C_1^2 K_1$.
Ainsi $\ln r_{i} - \ln r_{i-1} \leq - K_2 \on{aire} T(w_i)$.
Donc, en notant $r$ le rayon conforme du revêtement universel de $U$,
\[\ln r \leq - K_2 \sum_{i \in I} \on{aire} T(w_i))\]
et donc d'après le lemme~\ref{lem_vitalmajor},
\[\ln r \leq - \frac{K_2}{C_3} \on{aire} \bigcup_{e^w \in\bbD\setminus U}
T(w) = - K a\]  

La constante $K = 2\pi K_2/C_3$ que nous avons obtenue est très petite
(elle vaut approximativement $1/20000$), mais nous n'avons pas cherché à
optimiser le calcul. 


\bigskip\bigskip

\textsl{\large Minoration de $r$}

\bigskip

Cette minoration n'interviendra pas dans la correction de la thèse.

\

Le lemme suivant peut avoir son utilité dans d'autres contextes~:

\begin{lem}\label{lem_cvx}
 Soit $W \subsetneq \bbC$ un disque topologique contenant $0$.
 Pour tout fermé $A$ ne contenant pas $0$ tel que la composante
 connexe $U$ de $W \setminus A$ contenant $0$ est simplement connexe,
 soit $\lambda(A) = r(U)/r(W)$ où $r(\cdot)$ désigne le rayon conforme
 par rapport à $0$.
 \\
 Si $A$ et $B$ sont deux fermés vérifiant les conditions ci-dessus,
 alors $A\cup B$ les vérifie, et
 \[\lambda(A \cup B) \geq \lambda(A) \lambda(B)\]
\end{lem}
\begin{dem}
 Soit $\bbS^2$ la sphère de Riemann. Nous utiliserons les deux
 résultats de topologie suivants. Un ouvert connexe et simplement
 connexe de $\bbS^2$ a son complémentaire connexe. Toute composante
 connexe du complémentaire d'un fermé connexe de $\bbS^2$ est
 simplement connexe. 
 Soit maintenant $U$ associé à $A$ et $U'$ à
 $B$. Par hypothèse, $\bbC \setminus U$ est connexe et $\bbC \setminus
 U'$ aussi. Comme ils contiennent tout deux $\bbC\setminus W$, leur
 réunion $\bbC \setminus (U \cap U')$ est connexe. La composante
 connexe $U''$ de $W \setminus (A \cup B)$ qui contient $0$ est
 également la composante connexe de $U \cap U'$ qui contient
 $0$. Donc elle est simplement connexe.

 Considérons la fonction de Green $g$ de $W$ en $0$, et $g_A$
 (resp. $g_B$, $g_{A \cup B}$) celle de la composante de $W \setminus
 A$ contenant $0$ en $0$ (et ainsi de suite).
 Par fonction de Green nous entendons le supremum des fonctions
 sous-harmoniques positives s'annulant au bord et inférieures à
 $\on{cste} -\ln|z|$ en $0$~: c'est une fonction harmonique positive
 s'annulant au bord et ayant en $0$ pour développement limité
 $-\ln(z)+\ln(r)+o(1)$ où $r$ est le rayon conforme.
 Soit la fonction $h = g + g_{A \cup B} - g_A - g_B$, définie sur
 $U''$. Cette fonction
 est harmonique, même en $0$. L'inégalité de l'énoncé équivaut à $h(0)
 \geq 0$ (notons que nous avons, d'après les inclusions entre les divers
 ensembles, $g \geq g_A \geq g_{A \cup B}$ et $g \geq g_B \geq g_{A
 \cup B}$, mais que cela n'implique pas $h(0) \geq 0$).
 Par le principe du minimum, il suffit de montrer que les
 valeurs d'adhérences de $h$ au bord de $U''$ sont $\geq 0$.
 Tout point $z$ au bord de $U''$ est dans le bord d'au moins un des trois
 ensembles $W$, $U$, ou $U'$. Au bord de $W$, les quatre fonctions
 tendent vers $0$. Au bord de $U$, en termes de valeurs d'adhérences,
 $g_A$ est nulle, $g \geq g_B$ et $g_{A\cup B} \geq 0$, donc $h \geq
 0$. De même pour le bord de $U'$.
\end{dem}

On peut formuler ce résultat ainsi~: pour les domaines simplement
connexes, quand on enlève un même morceau à deux domaines l'un
contenant l'autre, la perte relative de rayon conforme est moindre
pour le domaine le plus petit.

\

Reprenons la suite $z_i = \exp(w_i)$ où $w_i$ est fournie par le
lemme~\ref{lem_vitally_zero}, et soit $V_5(z_i) = \exp(T_5(z_i))$.
L'ouvert $\ds U'=\bbD \setminus \bigcup_{i \in I} V_5(z_i)$ est
inclus dans $U$. Son rayon conforme est donc inférieur à celui de
$U$. D'autre part, il est simplement connexe (car étoilé par rapport à
$0$).
D'après le lemme~\ref{lem_cvx}, son rayon conforme est supérieur au
produit des pertes induites par chaque $V_5(z_i)$ pris isolément.
Comme dans le lemme~\ref{lem_calculs}, la contribution d'un seul
triangle est de l'ordre de son aire~:

\begin{lem}\label{lem_geodortho}
  Soit $u>0$, $x=\exp(-u)$, $U$ la composante connexe contenant $0$ de
  $\bbD$ privé de la géodésique passant par $x$ et orthogonale à
  $(0x)$, et $r$ le rayon conforme associé à $U$. Alors
  \[ r = 1 / \cosh(u) \]
  Quand $u \tend 0$, 
  \[ r = 1 - \frac{u^2}{2} + O(u^3) \]
\end{lem}

\begin{figure}[htbp]
 \begin{center}
  \begin{picture}(124,130)
   \scalebox{0.333}{
    \put(0,0){\includegraphics[120pt,200pt][490pt,590pt]{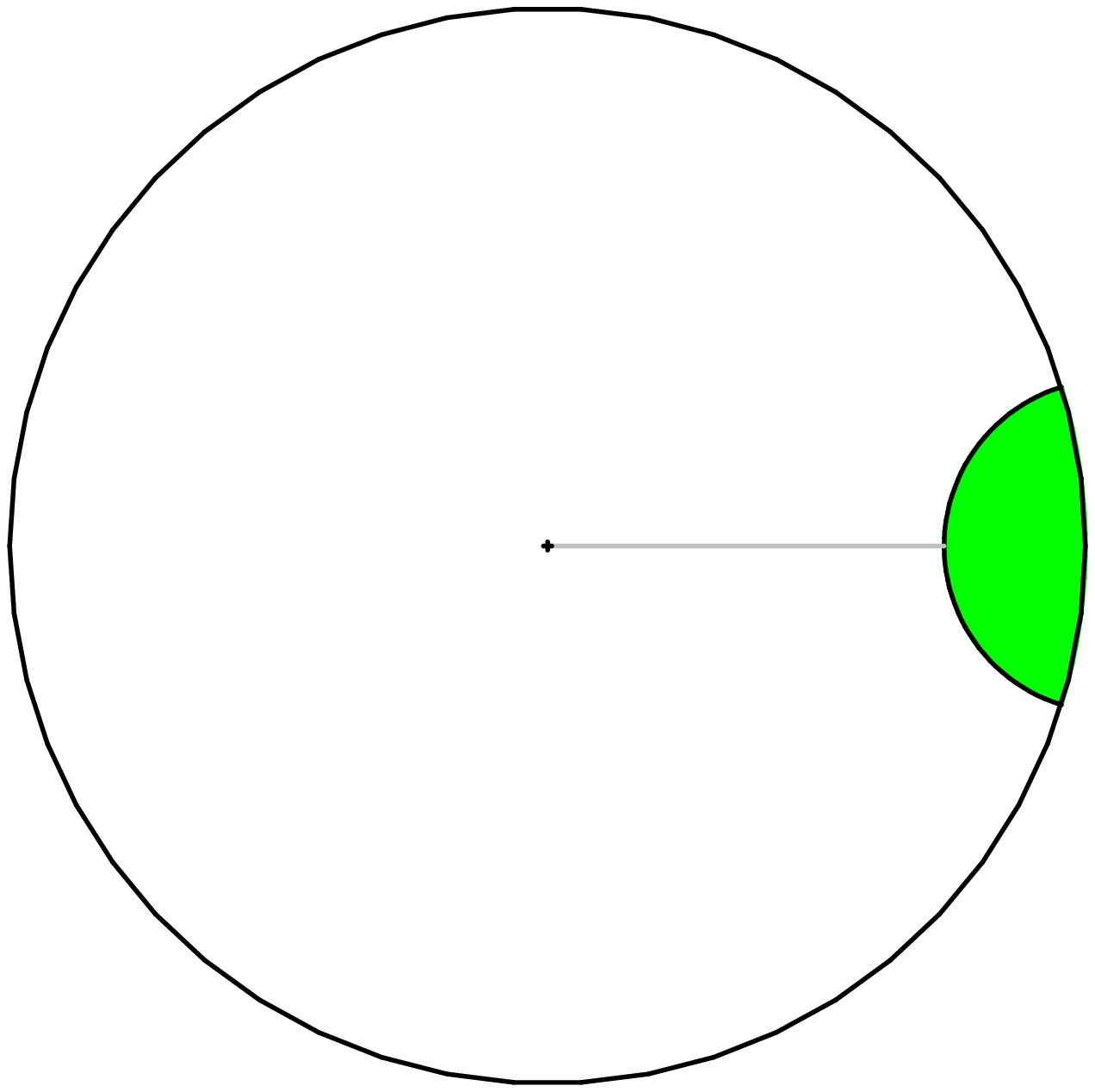}}
   }
  \end{picture}
 \end{center}
 \label{fig_geod_deux}
\end{figure}

\begin{lem}[contribution d'un triangle]\label{lem_vcinq}
  Il existe des constantes $K_3> K_4>1$ telle que
  si $0<u\leq\pi$, $x=e^{-u}$, $U=\bbD\setminus\{V_5(x)\}$,
  et $r$ est le rayon conforme associé à $U$, alors
  \begin{floatingfigure}[r]{60pt}
   \begin{picture}(60,117)
    \put(15,0){\scalebox{0.3}{
     \includegraphics[290,200][440,590]{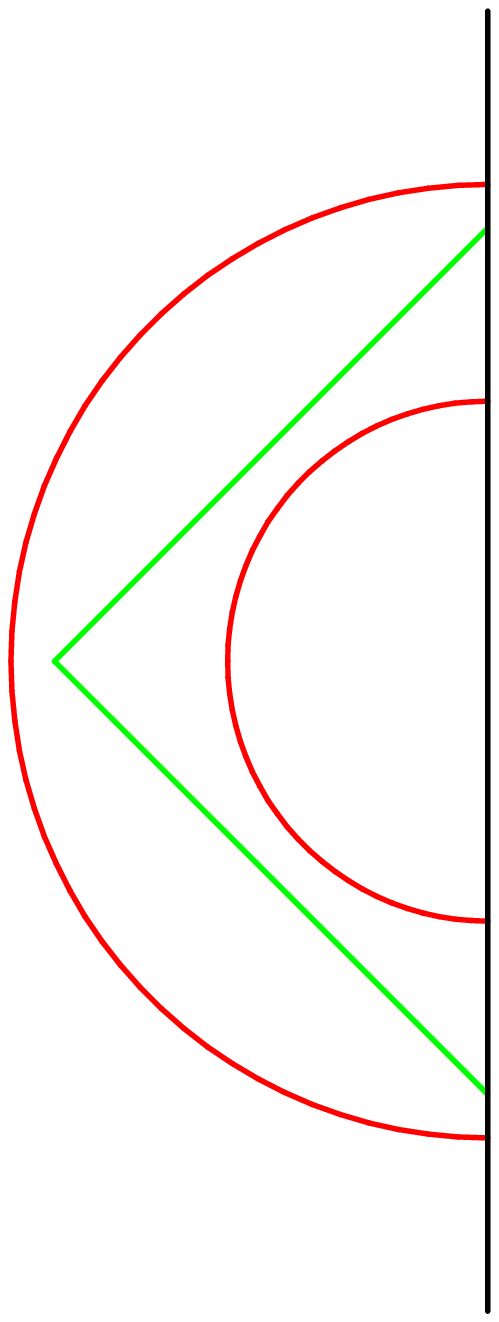} 
    }}
    \end{picture}
   \label{fig_geods}
  \end{floatingfigure}
  \[ -\ln(r)/\on{aire}\ T_5(u) \in [K_4,K_3]\]
\end{lem}
\begin{dem}
  Par un argument de continuité, il suffit de s'intéresser au cas où
  $u$ est proche de $0$. On peut alors placer le bord de $V_5(x)$ qui
  est très proche d'un triangle entre deux géodésiques comme sur la
  figure ci-contre. Avec des bases commensurables, leur aire est
  commensurable. On conclut ensuite avec le lemme~\ref{lem_geodortho}.
  Nous trouvons que $K_3 = 8/2\pi$ convient.
\end{dem}

D'après le lemme~\ref{lem_cvx} 
\[\ln r' \geq \sum_{i \in I} - K_3 \on{aire} T_5(w_i)\]
d'où
\begin{eqnarray*}
 \ln r \quad \geq \quad \ln r'
           & \geq & - K_3 \sum_{i \in I} \on{aire} T_5(w_i) \\
           & \geq & - K_3 C_3 . \on{aire}
             \bigcup_{\exp(w) \in \bbD\setminus U} T(w) \\
\end{eqnarray*}
d'après le lemme~\ref{lem_vitalminor}.

\

On obtient ainsi une constante $K' = 2\pi K_3 C_3$. Dans nos
estimation, $K'=200$, ce qui n'est pas excellent.



\clearpage 

\textbf{\LARGE III. Application à la correction de la thèse}

\bigskip

\medskip

Nous n'utiliserons que la majoration de $r$ dans la
proposition~\ref{prop_CRZ}.

\

\begin{lem}\label{lem_cone}
  Pour tout $u>0$, soit $T$ le triangle plein fermé de sommets
  $-u$, $\ic \pi$ et $-\ic \pi$. La plus petite
  longueur hyperbolique $l(u)$ dans $\bbD$ d'un chemin évitant $0$ et
  dont le relevé relie $T$ à un translaté de $T+2\ic\pi$, tend vers
  $+\infty$ quand $u \tend 0$.
\end{lem}

Nous le reformulerons ainsi~: pour passer le cône de paramètre $u$, un
chemin doit avoir une longueur $\geq l(u)$.

\begin{lem}
  Pour tout $d>0$, il existe une constante $K(d) > 1$ telle que pour
  tout ensemble fini de points $z_0$,$z_1$, \ldots, $z_n=z_0$ de
  $\bbD$, non nuls, tels que $z_{i+1}$ et $z_i$ sont reliés par un
  chemin $\gamma_i$ évitant $0$, de longueur hyperbolique $\leq d$
  dans $\bbD$, et tels que la concaténation des $\gamma_i$ fait au
  moins un tour autour de $0$, si on note $\mcZ = \{z_i\}_{i=1\ldots
  n}$ et $r'$ le rayon conforme de $\bbD \setminus \mcZ$, alors
  \[\ln r'\ \leq\ \frac{ \ln(\max |z_i| )}{K(d)} \ (<0)\]
\end{lem}
\begin{dem}
  Autrement dit, si on note $m= \max |z_i|$,
  \[ \frac{-\ln r'}{-\ln m }  \geq
  \frac{1}{K(d)}\]
  Premier cas~: l'un des points de $\mcZ$ est à distance $<
  e^{-\pi}$ de $0$. Soit $z_k$ le point le plus proche de $0$. Notons
  $|z_k| =\exp(-u)$~: $\ln m \geq -u$.
  D'autre part, $r'$ est plus petit que le rayon conforme de $\bbD
  \setminus \{z_k\}$, donc $\ln r' < \ln(u / \sinh u)$, d'où
  $-\ln r'/-\ln m > \ln(\sinh u/u) /u >
  \ln(\sinh(\pi)/\pi)/\pi > 0.4$, car $u > \pi$ et $u \mapsto
  \ln(\sinh u/u) /u$ est une fonction strictement croissante.
  \\
  Deuxième cas~: tous les points de $\mcZ$ sont à distance $\geq
  e^{-\pi}$. Soit $u>0$ tel que $u<\pi$ et $l(u)>d$ dans le
  lemme~\ref{lem_cone}.  A un élément $z \in \mcZ$, $z = \exp(w)$, on
  associe un intervalle de $\ic\bbR / \ic 2\pi \bbZ$, de centre $\ic
  \Im(w)$ et de longueur $2\pi |\Re(w)|/u$, noté $I(z)$. La réunion
  des $I(z)$ pour $z\in \mcZ$ recouvre tout le quotient $\ic\bbR / \ic
  2\pi \bbZ$.  En effet, supposons par l'absurde qu'un point $\ic y$
  échappe à cette réunion. Chaque $\gamma_i$ a une longueur
  hyperbolique dans $\bbD$ qui est $\leq d$, donc $< l(u)$. Or la
  succession des $\gamma_i$ fait au moins un tour autour de $0$, donc
  doit passer d'un coté à l'autre du cône issu de $\ic y$ et de
  paramètre $u$, or d'après le lemme précédent il est trop court pour
  cela.
  \\
  On sélectionne maintenant un sous-ensemble $E$ de $\mcZ$ dont les
  intervalles $I(z)$ sont deux à deux disjoints, en commençant par le
  plus grand, puis en sélectionnant par récurrence le plus grand parmi
  ceux qui restent et qui sont disjoints de ceux déjà sélectionnés.
  Alors, la réunion des intervalles de même centres et de longueur
  triple recouvre tout $\ic\bbR / \ic 2\pi \bbZ$.  Comme les éléments
  $I(z)$ pour $z \in E$ sont deux à deux disjoints, et contiennent la
  base de $T(z)$, les $T(\ln z)$ pour $z\in E$ sont deux à deux
  disjoints. Donc l'aire de leur réunion est égale à la somme de leurs
  aires, elle même égale à
 \begin{eqnarray*}
  \sum_{z\in E} |\Re(w)|^2
    & = & \sum_{z\in E} \frac{u |I(z)|}{2\pi} |\Re(w)| \\
    & \geq & \sum_{z\in E} \frac{u |I(z)|}{2\pi} |\ln m|\\
    & \geq & \frac{u |\ln m|}{2\pi} \sum_{z\in E} |I(z)| \\
    & \geq & \frac{u |\ln m|}{2\pi} \frac{2\pi}{3} \\
    & \geq & \frac{u}{3} |\ln m| \\
 \end{eqnarray*}
 Ainsi l'aire de $\ds
 \bigcup_{\exp(w) \in \mcZ} T(w)$ est minorée par une valeur
 proportionnelle à $|\ln m|$ (avec un coefficient dépendant de $d$).
 Il suffit alors d'appliquer la proposition~\ref{prop_CRZ}.
\end{dem}

\begin{lem}
  Pour tout $d>0$, pour toute fonction $f : \bbD \dans \bbD$ fixant
  $0$ et tous $r\in]e^{-1},1[$ et $q\in\bbN$, $q\geq 2$, tels que la
  distance hyperbolique dans $\bbD$ entre $r$ et $r \exp(\ic 2\pi/q)$
  est $\leq d$, alors
  \[\ln(r') \leq \ln(r)/K(d)\ (<0)\]
  où $r'$ est le rayon conforme de $\bbD \setminus f(r \bbU_q)$, et
  $K(d)$ est la constante du lemme précédent.
\end{lem}
\begin{dem}
  Il suffit de prendre pour $\gamma_i$ l'image par $f$ de la
  géodésique qui relie $r \exp(\ic 2\pi i/q)$ à $r \exp(\ic 2\pi
  (i+1)/q)$. Si l'un d'entre eux touche $0$ il suffit de
  considérer une valeur de $r$ légèrement inférieure et de passer à la
  limite. Le nombre de tours autour de $0$ de la concaténation de ces
  chemins est égal à au nombre de fois que $f$ atteint $0$ et donc est
  $>0$.
\end{dem}

Pour conclure, remarquons que pour toute valeur de $d$, il existe une
constante $R(d)>0$ telle que si $q\geq 3$, $r \geq e^{-1}$, (valeurs
que j'ai fixées arbitrairement), et si $q|\ln(r)| \geq R(d)$, alors la
condition sur les distances hyperboliques des lemmes précédents est
vérifiée.

\clearpage

\textbf{\LARGE Conclusion}

\bigskip

\medskip

Je pense que les constantes $K'$ et $K$ peuvent être considérablement
améliorées. On peut également chercher à remplacer la forme
triangulaire par une autre, pour avoir un quotient $K'/K$ encore
meilleur (dans le présent article, nous obtenons un quotient de
l'ordre de $4.10^6$).

\bigskip

\bigskip

\textbf{\LARGE Remerciements}

\bigskip

\medskip

Je tiens à remercier Xavier~Buff, Adrien~Douady, et Michel~Zinsmeister
pour d'utiles discussions. La question de départ, pendant que j'en
rédigeait une preuve (erronée), a inspiré à Douady
une question différente analogue. Rohde et Zinsmeister l'ont résolue.
En présentant ma démonstration à Buff, nous nous sommes rendu compte
d'une erreur dans la preuve d'un des lemmes. 
Le présent addendum reprend certaines des idées de Rohde et
Zinsmeister pour corriger cela. En même temps, Buff a exploré une
approche plus élégante pour corriger ce lemme, qu'il a menée à bien en
démontrant un résultat général fort élégant.

\end{document}